\documentclass[10pt,preprint]{amsart}
\usepackage{amssymb,amsfonts,amsthm,amsmath,amscd,relsize}
\usepackage{hyperref}
\usepackage{enumerate}
\usepackage{comment}

 \usepackage{epsfig}
 \usepackage{graphicx}
 \usepackage{lipsum}
\usepackage{xcolor}
\usepackage{float}
\usepackage{soul}
\usepackage{adjustbox}

\usepackage{enumitem}

\usepackage{pgf,tikz,pgfplots}
\pgfplotsset{compat=1.13}
\usepackage{mathrsfs}
\usetikzlibrary{arrows}

\definecolor{qqqqff}{rgb}{0,0,1}
\definecolor{qqwuqq}{rgb}{0,0.39215686274509803,0}
\definecolor{uuuuuu}{rgb}{0.266,0.266,0.266}

\def\k{\kappa}

\def\A{\textbf{\textit{A}}}
\def\B{\textbf{\textit{B}}}
\def\C{\textbf{\textit{C}}}
\def\D{\textbf{\textit{D}}}
\def\X{\textbf{\textit{X}}}
\def\k{\textbf{\textit{k}}}

\theoremstyle{plain}
\newtheorem*{theorem}{Theorem}

\theoremstyle{definition}

\theoremstyle{remark}

\makeatletter
\@addtoreset{footnote}{page}
\makeatother

\definecolor{ttzzqq}{rgb}{0.2,0.6,0}

\begin{document}

\title{A vector identity for quadrilaterals}

\author{Christian Aebi and Grant Cairns}

\address{Coll\`ege Calvin, Geneva, Switzerland 1211}
\email{christian.aebi@edu.ge.ch}
\address{Department of Mathematics, La Trobe University, Melbourne, Australia 3086}
\email{G.Cairns@latrobe.edu.au}

%\begin{abstract} We show ...
%\end{abstract}

\maketitle

Consider a quadrilateral $ABCD$ in the plane. We let $K_{ABCD}$ denote its area and we use obvious similar notation for the area of various triangles. Consider the position vectors $\A,\B,\C,\D$ of the vertices. Then one has the following result, which one might call the \emph{Jacobi identity for quadrilaterals}.

\begin{theorem}\ 
$K_{BCD}\, \A -K_{ACD}\, \B+K_{ABD}\, \C  -K_{ABC}\, \D=0$.
\end{theorem}

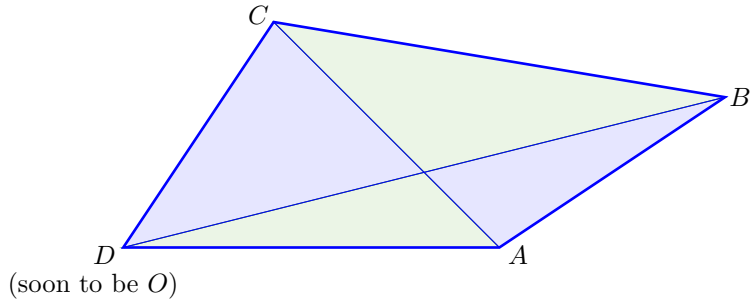
\begin{center}
\begin{tikzpicture}[scale=.5][line cap=round,line join=round,>=triangle 45,x=1cm,y=1cm]
\draw[color=ttzzqq,fill=ttzzqq,fill opacity=0.1] (0,0) -- (10,0)-- (8,2)-- cycle; 
\draw[color=ttzzqq,fill=ttzzqq,fill opacity=0.1] (8,2) -- (4,6)-- (16,4)-- cycle; 
\draw[color=qqqqff,fill=qqqqff,fill opacity=0.1] (8,2) --  (10,0) -- (16,4)   -- cycle; 
\draw[color=qqqqff,fill=qqqqff,fill opacity=0.1] (8,2) --  (0,0) -- (4,6)   -- cycle; 
\draw[line width=1pt,color=qqqqff] (0,0) -- (10,0)--(16,4)--(4,6)-- cycle; 
%\draw[line width=1pt,color=qqqqff]   (12,0)-- (12,3);
%\draw[line width=1pt,color=qqqqff]   (12,3)-- (72/7,30/7);
\draw[color=black] (-.5,-.2) node {$D$};
\draw[color=black] (-.8,-1.) node {(soon to be $O$)};
\draw[color=black] (10.5,-.2) node {$A$};
\draw[color=black] (16.4,4) node {$B$};
\draw[color=black] (3.6,6.2) node {$C$};
\end{tikzpicture}
%\caption{A quadrilateral}\label{F:quad}
\end{center}

\begin{proof} First notice that the identity is unchanged by translation. Indeed, for a translation by a vector $\X$, the areas are unchanged and
\begin{align*}
K_{BCD}\, (\A+\X) &-K_{ACD}\, (\B+\X)+K_{ABD}\, (\C+\X)  -K_{ABC}\, (\D+\X)\\
%&=K_{BCD}\, \A -K_{ACD}\, \B+K_{ABD}\, \C  -K_{ABC}\, \D +(K_{BCD}-K_{ACD}+K_{ABD}  -K_{ABC} )+\X\\.
&=K_{BCD}\, \A -K_{ACD}\, \B+K_{ABD}\, \C  -K_{ABC}\, \D,
\end{align*}
since $K_{BCD}+K_{ABD} =K_{ABCD}=K_{ACD} +K_{ABC}$. So we may take $D$ to be the origin, $O$. Then our desired identity is
\[
K_{BCO}\, \A -K_{ACO}\, \B +K_{ABO}\, \C =0.
\]
Observe that the vector cross product $\B\times \C$ has magnitude $2K_{BCO}$ and direction $\k$ perpendicular to the plane. Further, $\k\times \A$ is the vector obtained from $\A$ by rotating it through angle $\frac{\pi}2$ in the positive direction; let us denote this vector $\A'$. Similarly, we define $\B'$ and $\C'$. Now, rotating the vector $K_{BCO}\, \A -K_{ACO}\, \B +K_{ABO}\, \C$ through $\frac{\pi}2$, and multiplying by 2,  we obtain
\begin{align*}
2(K_{BCO}\, \A' &-K_{ACO}\, \B' +K_{ABO}\, \C')\\
&= (\B\times \C)\times \A- (\A\times \C)\times \B+(\A\times \B)\times \C=0,
\end{align*}
by the Jacobi vector triple product identity. Hence $K_{BCO}\, \A -K_{ACO}\, \B +K_{ABO}\, \C =0$, as desired.
\end{proof}

It remains for us to clarify what we mean by ``quadrilateral''. In fact, the word does not have an entirely universal meaning; see \cite{UG}. However, the astute reader will have noticed that the above proof doesn't require any property of $ABCD$; it is actually true for every choice of 4 points $A,B,C,D$ in the plane, provided one interprets the numbers $K_{ABC}$ etc as \emph{signed areas}; that is,   $K_{ABC}$ is positive if $A,B,C$ are positively oriented, and negative otherwise.

For more on the wonderful world of quadrilaterals, see the delightful book by Alsina  and  Nelsen \cite{AN}.

%%%%%%%%%%%%%%%%%%%%%%%%%%%%%%%%%%%%%%%%%%%%%%%%%%%%
\bibliographystyle{amsplain}
{}

%%%%%%%%%%%%%%%%%%%%%%%%%%%

\end{document}